\documentclass[reqno,10pt]{amsart}
\usepackage[english]{babel}
\usepackage{amsmath}
\usepackage{amssymb}
\usepackage[latin1]{inputenc}
\usepackage{dsfont}
\usepackage[bookmarksnumbered,colorlinks]{hyperref}
\usepackage{graphics}
\usepackage{enumerate}
\usepackage[all]{xy}
\usepackage{mathrsfs}

\newcommand{\R}{\mathds R}
\newcommand{\N}{\mathds N}
\newcommand{\Z}{\mathds Z}

\newcommand{\C}{{\mathcal C}}
\newcommand{\Cb}{{\mathfrak C}}

\begin{document}

\newtheorem{theorem}{Theorem}[section]
\newtheorem{proposition}[theorem]{Proposition}
\newtheorem{lemma}[theorem]{Lemma}
\newtheorem{corollary}[theorem]{Corollary}
\theoremstyle{definition}
\newtheorem{definition}[theorem]{Definition}
\newtheorem{remark}[theorem]{Remark}


%
%

\title{Anisotropic tensor calculus}

\author{MIGUEL \'ANGEL JAVALOYES}
\address{Department of Mathematics, University of Murcia\\
Facultad de Matem\'aticas, Campus de Espinardo, Murcia, 30100, Spain\\
majava@um.es}


\maketitle


\begin{abstract}
We introduce the anisotropic tensor calculus, which is a way of handling with tensors that depend on the direction remaining always in the same class. This means that the derivative of an anisotropic tensor is a tensor of the same type. As an application we show how to define derivations using anisotropic linear connections in a manifold. 
In particular, we show that the Chern connection of a Finsler metric can be interpreted as the Levi-Civita connection and we introduce the anisotropic curvature tensor. We also relate the concept of anisotropic connection with the classical concept of linear connections in the vertical bundle.  Furthermore, we also introduce the concept of anisotropic Lie derivative.
\end{abstract}

\section{Introduction}

Anisotropic tensors appear in situations in that there is a strong dependence on the direction as it is the case of Finsler metrics and sprays. Traditionally, derivations of these tensors in a manifold $M$ have been handled with linear connections in the vertical subbundle of $TTM$  over the slit tangent bundle $TM\setminus 0$. This has the following drawback: you need to carry some information about derivatives which is not geometrically relevant. The most apparent consequence of this approach is that there are many linear connections which can be associated with a Finsler metric, losing the unicity of the Levi-Civita connection in Riemannian geometry. Our aim is to solve these limitations by considering an anisotropic linear connection $\nabla$, which is a connection in the manifold but with a dependence on the direction (see Definition \ref{aniconnection}). Then when we compute the derivative of an anisotropic tensor $T$ in $v\in TM\setminus 0$, we choose an extension $V\in \mathfrak{X}(M)$ of $v$ and make the computation with the affine linear connection $\nabla^V$ and the tensor $T_V$ (fixing the direction $V(p)$ at every point  $p$ of the manifold $M$). In principle, the result depends on the extension $V$, but this dependence disappears when a suitable vertical derivative of the tensor is substracted (see \eqref{formderi} and \eqref{partialTexpres}). This procedure has been inspired by the interpretation of the Chern connection as a family of affine connections, which was first developed in \cite{Mat80}. Later, in \cite{Jav14a,Ja14b}, the relation of the flag curvature of the Finsler metric with the family of affine connections was completely clarified, setting down the foundations of anisotropic geometry (see also \cite{JaSo14}). The central result of this paper is Theorem \ref{existderiv}, which allows us to deal with and to develop the anisotropic tensor calculus departing from an anisotropic connection (which is a concept known in literature \cite[Chapter 7]{Sh01}).  The other central concept is the curvature tensor of an anisotropic connection introduced in \S  \ref{curvanis}. 

One of the main virtues of this approach is that the treatment of Finsler Geometry parallels the classical one of Riemannian Geometry using the Levi-Civita connection.
In fact, with our anisotropic tensor calculus, the problem of unicity for connections associated with Finsler metrics and sprays is completely solved. As a matter of fact, given a linear connetion on the vertical subbundle, we can construct an anisotropic linear connection on the manifold (see \S \ref{relationbundles}). It turns out that the four classical linear connections associated with a Finsler metric: Berwald, Cartan, Chern and Hashiguchi, project into only two anisotropic connections. The Berwald and Hashiguchi connections project into a connection that  is the most natural one to be associated with a spray, whereas the Chern and Cartan ones project into an anisotropic connection which is the Levi-Civita connection of a Finsler metric, since it is the unique anisotropic connection which is torsion-free and such that $\nabla g=0$, namely, the derivative of the fundamental tensor is trivial (see \S \ref{chernconnection}). Another remarkable fact is that we do not need to fix the nonlinear connection in order to obtain the invariants associated with a Finsler metric or a spray. All the invariants are obtained directly with the anisotropic connection differently from what happens when you consider a connection in the vertical subbundle, where the nonlinear connection is required in order to compute the geometric invariants.


In \S \ref{relationbundles}, it is clarified how the computations with anisotropic connections can be related with those of classical linear connections of the vertical bundle. In turns out that when an appropriate non-linear connection is fixed, a classical linear connections is determined by an anisotropic connection and a $(0,2)$-anisotropic tensor $\mathfrak C$, which takes control of the vertical quantities. Proposition \ref{relationcurder} provides the relation between the curvatures and derivatives of the linear connection in the vertical bundle and the derivatives and curvatures of the associated anisotropic connection and derivatives of the anisotropic tensor $\mathfrak C$. This concludes that with the approach of anisotropic connections no information is lost.

Our last achievement is to introduce the anisotropic Lie derivative. This theory has been studied in anisotropic geometry for a long time \cite{Yano57}, but without a clear tensorial development as in Riemannian geometry. The anisotropic approach allows one to define the Lie derivative of any anisotropic object by using the Lie bracket of vector fields as an  anisotropic derivation (see Definition \ref{ander}).  There is another approach that uses $TTM$ to define Lie derivatives, which entails more computations (see for example \cite{Lovas04}). In Proposition \ref{Lieflow}, we obtain a very natural interpretation of the Lie derivative in terms of the flow of a vector field.

The paper is organized as follows. In Section \ref{antensor}, we introduce the anisotropic tensor calculus,  the main result being Theorem \ref{existderiv}, which ensures the existence of a unique anisotropic tensor derivation for every anisotropic derivation as in Definition~\ref{ander}.  In \S \ref{anlinconn}, we give the definition of anisotropic linear connection and the associated curvature tensor. In \S \ref{nonlinearity} and \S \ref{curvanis}, we introduce the non-linear and Ehresmann connections and the curvature tensor associated with an anisotropic connection, concluding the section   with the introduction of the covariant derivative along a curve associated with an anisotropic connection in \S \ref{covariant}. 
In \S \ref{chernconnection}, we show that the Chern connection is the Levi-Civita connection of a Finsler metric, namely, it is the only anisotropic connection which is torsion-free and parallel. 
Furthermore,  in \S \ref{nonlinear} we introduce the notion of spray and its associated nonlinear connection. Some classical tensors are also introduced in \S \ref{berwaldcon} to describe the difference tensor between the Chern and Berwald connections. 
In \S \ref{relationbundles} it is shown how to associate an anisotropic connection to a connection in the vertical subbundle and how one can relate their curvature tensors. Finally, in \S \ref{lieder}, we show how to define the anisotropic Lie derivative using the Lie bracket and Theorem \ref{existderiv}. We conclude with an interpretation of the anisotropic Lie derivative using the flow of a vector field (see Proposition \ref{Lieflow}).

\section{Anisotropic Tensor Calculus}\label{antensor}
\subsection{Anisotropic tensors} 
Let $M$ be a manifold of dimension $n$, $TM$ its tangent bundle, $TM^*$ its cotangent bundle and $A\subset TM\setminus 0$, an open subset such that $\pi(A)=M$. We will use the notations $\pi:TM\rightarrow M$ and $\tilde \pi:TM^*\rightarrow M$,  for the natural projections of the tangent and cotangent bundle, respectively. Using the restriction  $\pi|_A:A\rightarrow M$, we can lift the vector bundle $\pi:TM\rightarrow M$ to $A$ obtaining a vector bundle $\pi_A^*(M)$ and by lifting $\tilde\pi:TM^*\rightarrow M$, we get another vector bundle $\tilde\pi_A^*(M)$:
\begin{displaymath}
    \xymatrix{ \pi_A^*(M) \ar[d]_{\pi^*_A} & TM \ar[d]^{\pi} \\
               A\subset TM\setminus 0 \ar[r]^{\quad\,\,\pi|_A}  & M  }
   \quad \quad \xymatrix{ \tilde\pi_A^*(M) \ar[d]_{\tilde\pi^*_A} & TM^* \ar[d]^{\tilde\pi} \\
               A\subset TM\setminus 0 \ar[r]^{\quad\,\,\pi|_A}  & M  }
\end{displaymath}
Observe that at $v\in A$, we have that $(\pi_A^*)^{-1}(v)= T_{\pi(v)}M$, namely, the fiber of $\pi_A^*(M)$ which projects to $v\in A$ is $T_{\pi(v)}M$, and also $(\tilde\pi_A^*)^{-1}(v)= T_{\pi(v)}M^*$. As a consequence, the dimension of $\pi_A^*(M)$ (and also of $\tilde\pi_A^*(M)$) as a manifold is $3n$.
 We can also make tensor products of these vector bundles as \[\pi_A^*(M)^{\otimes r}\otimes (\tilde\pi_A^*(M))^{\otimes s}:=\pi_A^*(M)\otimes \overbrace{\cdots}^{r} \otimes \pi_A^*(M)\otimes\tilde\pi_A^*(M)\otimes \overbrace{\cdots}^{s} \otimes \tilde\pi_A^*(M),\]
 with projections $(\pi_A^*)^r_s:\pi_A^*(M)^{\otimes r}\otimes (\tilde\pi_A^*(M))^{\otimes s}\rightarrow A$. Moreover, if $v\in A$, then $((\pi_A^*)^r_s)^{-1}(v)=(T_{\pi(v)}M)^{\otimes r}\otimes (T_{\pi(v)}M^*)^{\otimes s}$.
 \begin{definition}
Given non-negative integer numbers $r,s\in\N^*$ such that $r+s>0$, we define an $A$-anisotropic $(r,s)$-tensor in $M$ as a smooth section of $\pi_A^*(M)^{\otimes r}\otimes \tilde\pi_A^*(M)^{\otimes s}$, namely, a smooth map $T: A\rightarrow \pi_A^*(M)^{\otimes r}\otimes \tilde\pi_A^*(M)^{\otimes s}$ such that $(\pi^*_A)^r_s\circ T$ is the identity. 
\end{definition}
The space of all the $A$-anisotropic $(r,s)$-tensors of $M$ will be denoted by $\mathfrak{T}^r_s(M,A)$. By definition, $\mathfrak{T}^0_0(M,A)$ will be the space of smooth real functions on $A$,  $f:A\rightarrow \R$, also denoted as ${\mathcal F}(A)$. Given an $A$-anisotropic tensor $T\in \mathfrak{T}^r_s(M,A)$ and $v\in A$, we will use the notation $T_v=T(v)\in (T_{\pi(v)}M)^{\otimes r}\otimes (T_{\pi(v)}M^*)^{\otimes s}$. Given a system of coordinates $(\Omega,\varphi)$ in $M$ with $\varphi=(x^1,x^2,\ldots,x^n):\Omega\rightarrow U\subset  \R^n$, if we denote by $\partial_i$ the partial vector fields in $\Omega$ and $dx^i$ their dual one-forms, then for every $T\in \mathfrak{T}^r_s(M,A)$ and every $v\in A\cap T\Omega$, we have that
\[T_v=\sum_{\text{\tiny$\begin{array}{cc} 1\leq i_1,i_2,\ldots,i_r\leq n\\1\leq j_1,j_2,\ldots,j_s\leq n\end{array}$}} T_{j_1j_2\ldots j_s}^{i_1i_2\ldots i_r}(v)\, \left(\partial_{i_1}\otimes\cdots\otimes \partial_{i_r}\otimes dx^{j_1}\otimes\cdots \otimes dx^{j_s}\right)|_{\pi(v)},
\]
where $ T_{j_1j_2\ldots j_s}^{i_1i_2\ldots i_r}:A\cap T\Omega\rightarrow \R$.
 Observe that at every point $p\in M$, we can identify $T_pM$ with the dual vector space of $T_pM^*$ and $T_pM^*$, with the dual vector space of $T_pM$.  Then the elements in the fibers of $\pi_A^*(M)^{\otimes r}\otimes \tilde\pi_A^*(M)^{\otimes  s}$ can be identified with multilinear maps
\[\phi:(T_pM^*)^r\times (T_pM)^s\rightarrow \R.\]
 With this interpretation, $T_{j_1j_2\ldots j_s}^{i_1i_2\ldots i_r}(v)=T_v(dx^{i_1},dx^{i_2},\ldots,dx^{i_r},\partial_{j_1},\partial_{j_2}\ldots,\partial_{j_s})$. As in the classical case of tensors (see \cite[Proposition 2.2]{Oneill}), we can also identify $T\in \mathfrak{T}^r_s(M,A)$ with the ${\mathcal F}(A)$-multilinear map
 \begin{equation}\label{tensormultilinear}
 T:\mathfrak{T}^0_1(M,A)^r\times \mathfrak{T}^1_0(M,A)^s\rightarrow
 {\mathcal F}(A),
 \end{equation}
 where ${\mathcal F}(A)=\{f:A\rightarrow \R: f\in C^\infty\}$. This is because for any 
 \[(\theta^1,\theta^2,\ldots,\theta^r,X_1,\ldots,X_s)\in \mathfrak{T}^0_1(M,A)^{r}\times \mathfrak{T}^1_0(M,A)^{s}\] and $v\in A$, the ${\mathcal F}(A)$-multilinearity implies that 
 \[T(\theta^1,\theta^2,\ldots,\theta^r,X_1,\ldots,X_s)(v)\]
 depends only on $T$ and $\theta^1(v),\theta^2(v),\ldots,\theta^r(v),X_1(v),\ldots,X_s(v)$, and not on their particular extensions $\theta^1,\theta^2,\ldots,\theta^r,X_1,\ldots,X_s$. Therefore this implies that we can define a map $T_v:(T_{\pi(v)}M^*)^r\times (T_{\pi(v)}M)^s\rightarrow \R$ for every $v\in A$ with smooth dependence on $v$. Moreover, from now on, we will use the following notation: $T_v(\theta^1,\theta^2,\ldots,\theta^r,X_1,\ldots,X_s)=T(\theta^1,\theta^2,\ldots,\theta^r,X_1,\ldots,X_s)(v)$, and we will often consider $T_v$ evaluated in local fields around $\pi(v)$, namely, defined in a neighborhood of $\pi(v)$, or even, in vectors in $T_{\pi(v)}M$.
 
 In the case we have an ${\mathcal F}(A)$-multilinear map 
 \begin{equation}\label{forthecurv}
 T: \mathfrak{T}^1_0(M,A)^s\rightarrow
 \mathfrak{T}^1_0(M,A),
\end{equation} 
  we define 
 $\bar{T}:\mathfrak{T}^0_1(M,A)\times
 \mathfrak{T}^1_0(M,A)^s\rightarrow {\mathcal F}(A)$ 
by  
\begin{equation}\label{interpr}
 \bar{T}(\theta,X_1,\ldots,X_s)=\theta(T(X_1,\ldots,X_s)),
 \end{equation}
 which is a $(1,s)$-tensor field. We shall consider $T$ to be a tensor field itself, using the formula above only when necessary.
 \begin{remark}\label{extendT}
 We will consider sometimes $\mathfrak{X}(M)$ as a subset of  $\mathfrak{T}^1_0(M,A)$, namely, a vector field $X\in\mathfrak{X}(M)$  can be lifted to $X^*\in\mathfrak{T}^1_0(M,A)$ by making $X^*(v)=X(\pi(v))$, and we will identify $X^*\equiv X$. Analogously, we will consider the space of one-forms on $M$, denoted $\mathfrak{X}^*(M)$, as a subset of  $\mathfrak{T}^1_0(M,A)$. Moreover,
  to define a tensor $T \in\mathfrak{T}^r_s(M,A)$ is equivalent to define an ${\mathcal F}(M)$-multilinear map \[T:\mathfrak{X}^*(M)^r\times
  \mathfrak{X}(M)^s\rightarrow {\mathcal F}(A).\]
  Then we can extend this map to $\mathfrak{T}^0_1(M,A)^r\times \mathfrak{T}^1_0(M,A)^s$ using a (local) frame on $M$ and ${\mathcal F}(A)$-multilinearity. It is straightforward to check that this extension does not depend on the frame. In a similar way, it is enough to define a map $T:\mathfrak{X}(M)^s\rightarrow  \mathfrak{T}^1_0(M,A)$, for determining a tensor as in \eqref{forthecurv}.
 \end{remark}	
 
 In the following, we will omit the subset $A$ when it is clear by the context. So we will speak just about anisotropic tensors.
Let us relate the notion of anisotropic tensor with the classical notion of tensor on a manifold. 
\begin{definition}\label{partialT}
We say that a vector field $V$ in an open subset $\Omega\subset M$ is $A$-admissible (with $A\subset TM\setminus 0$) if $V(p)\in A$ for every $p\in \Omega$.
\end{definition}
Observe that when we fix an $A$-admissible vector field $V$ on an open subset $\Omega\subset M$, then for every $T\in \mathfrak{T}^r_s(M,A)$ we can define a tensor $T_V\in \mathfrak{T}^r_s(\Omega)$ as
\[T_V:{\mathfrak X}^*(\Omega)^r\times {\mathfrak X}(\Omega)^s\rightarrow {\mathcal F}(\Omega)\]
such that $T_V(\theta^1,\ldots,\theta^r,X_1,\ldots,X_s)(p)=
T_{V(p)}(\theta^1(p),\ldots,\theta^r(p),X_1(p),\ldots,X_s(p))$, where ${\mathcal F}(\Omega)=\{f:\Omega\rightarrow \R: f\in C^\infty\}$. As usual, ${\mathfrak X}(\Omega)$ and ${\mathfrak X}^*(\Omega)$ denote, respectively, the space of vector fields and one-forms on $\Omega$.
On the other hand, the dependence on directions makes it necessary to define derivatives in the vertical bundle.
\begin{definition}\label{defi:vert}
Given an anisotropic tensor $T\in \mathfrak{T}^{r}_s(M,A)$, we define its {\em vertical derivative} as the tensor  $\partial^\nu T\in \mathfrak{T}^{r}_{s+1}(M,A)$ given by
 \[(\partial^\nu T)_v(\theta^1,\theta^2,\ldots,\theta^r,X_1,\ldots,X_s,Z)=\frac{\partial}{\partial t} T_{v+tZ(v)}(\theta^1,\theta^2,\ldots,\theta^r,X_1,\ldots,X_s)|_{t=0}\]
 for any $(\theta^1,\theta^2,\ldots,\theta^r,X_1,\ldots,X_s,Z)\in \mathfrak{X}^*(M)^{r}\times \mathfrak{X}(M)^{s+1}$ (recall Remark \ref{extendT}).
 \end{definition}
In the following, when working in the tangent bundle $TM$, we will use the natural coordinates $(T\Omega,\tilde{\varphi})$ associated with a coordinate system $(\Omega,\varphi)$ of $M$. As usual, we will denote the coordinates of a point $v\in T\Omega$ as 
\begin{equation}\label{naturalcoord}
\tilde\varphi=(x,y)=(x^1,x^2,\ldots,x^n,y^1,y^2,\ldots,y^n).
\end{equation}
We will also use the {\it  Einstein summation convention} when possible and will omit the coordinate functions $\varphi$ and $\tilde\varphi$ to avoid clutter in equations. Let us see how one can compute the coordinates of $\partial^\nu T$ in $(\Omega,\varphi)$:
\begin{equation}\label{coordvertical}
(\partial^\nu T)_{j_1\ldots j_sj_{s+1}}^{i_1\ldots i_r}=\frac{\partial T_{j_1\ldots j_s}^{i_1\ldots i_r}}{\partial y^{j_{s+1}}},
\end{equation} 
where the partial derivative is computed using the natural coordinates in $TM$ given in \eqref{naturalcoord}. Observe that as we consider the coordinates $ T_{j_1\ldots j_sj_{s+1}}^{i_1\ldots i_r}$ as functions on $T\Omega\cap A$, we are omitting some coordinate maps. In particular, we are using the following identification:
\[\frac{\partial T_{j_1\ldots j_s}^{i_1\ldots i_r}}{\partial y^{j_{s+1}}}=\frac{\partial( T_{j_1\ldots j_s}^{i_1\ldots i_r}\circ \tilde\varphi^{-1})}{\partial y^{j_{s+1}}}\circ \tilde\varphi,\]
which will be used from now on without further comments.

 We say that an anisotropic tensor $T\in \mathfrak{T}^{r}_s(M,A)$ is {\it positive homogeneous of degree $l\in\Z$} if $T_{\lambda v}=\lambda^l T_v$ for every $v\in A$ and $\lambda>0$. In such a case, it follows straightforwardly that
 \begin{equation}
 (\partial^\nu T)_{ v}(\theta^1,\theta^2,\ldots,\theta^r,X_1,\ldots,X_s,V)=l T_v(\theta^1,\theta^2,\ldots,\theta^r,X_1,\ldots,X_s),
 \end{equation}
where $V\in {\mathfrak X}(M)$ is an extension of $v$.
Like in the classical case of tensors on a manifold, we can define tensor products: if $T_1\in \mathfrak{T}^r_s(M,A)$ and $T_2\in \mathfrak{T}^{r'}_{s'}(M,A)$, then $T_1\otimes T_2 \in \mathfrak{T}^{r+r'}_{s+s'}(M,A)$ and it is given by
\begin{multline*}
(T_1\otimes T_2)(\theta^1,\ldots,\theta^{r+r'},X_1,\ldots,X_{s+s'})\\=
T_1(\theta^1,\ldots,\theta^r,X_1,\ldots,X_s)T_2(\theta^{r+1},\ldots,\theta^{r+r'},X_{s+1},\ldots,X_{s+s'}).
\end{multline*}
Observe that when one of the factors belongs to $ \mathfrak{T}^{0}_{0}(M,A)\equiv {\mathcal F}(A)$, the tensor product is just the usual product by a function.
 When both $r,s\geq 1$, we also can define contractions.
\begin{definition}
 Given an $(r,s)$-tensor $T$, its $(i,j)$-contraction is defined as the $(r-1,s-1)$-tensor ${\bf C^i_j}(T)$ such that in a system of coordinates $(\Omega,\varphi)$ as above, its coordinates are given by 
\[{\bf C^i_j}(T)_{j_1j_2\ldots j_{s-1}}^{i_1i_2\ldots i_{r-1}}(v)=\sum_{m=1}^n T_{j_1j_2\ldots j_{j-1}m j_{j}\ldots j_{s-1}}^{i_1i_2\ldots i_{i-1}m i_{i}\ldots i_{r-1}}(v),\]
for every $v\in T\Omega\cap A$.
It is not difficult to check that it is well-defined, namely, it does not depend on the coordinate system. Moreover, for a tensor $X\otimes \theta\in \mathfrak{T}^{1}_{1}(M,A)$, with $X\in  \mathfrak{T}^1_0(M,A)$ and $\theta \in  \mathfrak{T}^0_1(M,A)$,  contractions are defined as ${\bf C^1_1}(X\otimes \theta)=\theta(X)$.
\end{definition}
\subsection{Tensor derivations}
\begin{definition}
A tensor derivation $\mathcal D$ is an $\R$-linear map 
\[{\mathcal D}:\bigcup_{r\geq 0,s\geq 0}\mathfrak{T}^r_s(M,A)\rightarrow \bigcup_{r\geq 0,s\geq 0}\mathfrak{T}^r_s(M,A),\]
such that the restriction ${\mathcal D}^r_s={\mathcal D}|_{\mathfrak{T}^r_s(M,A)}$ has image included in $\mathfrak{T}^r_s(M,A)$
for any $r\geq 0, s\geq 0$, and it satisfies the following two properties:
\begin{enumerate}[(i)]
	\item  for any tensors $T_1$ and $T_2$:
${\mathcal D}(T_1\otimes T_2)={\mathcal D}T_1\otimes T_2+T_1\otimes {\mathcal D}T_2$,
\item   for any tensor $T$ and any contraction $\bf C$, ${\mathcal D}({\bf C}(T))={\bf C}({\mathcal D} (T))$.
\end{enumerate}
\end{definition}
It is easy to check that tensor derivations are local in the sense that if we consider an open subset $\mathcal A$ of $A$, then ${\mathcal D}_{\mathcal A}(T|_{\mathcal A})=({\mathcal D} T)|_{\mathcal A}$ and it obeys the product rule
\begin{align}\label{productrule}
{\mathcal D}(T(\theta^1,\ldots,\theta^r,X_1,\ldots,X_s))=& 
({\mathcal D}T)(\theta^1,\ldots,\theta^r,X_1,\ldots, X_s)\nonumber\\
&+\sum_{i=1}^r T(\theta^1,\ldots,{\mathcal D}\theta^i,\ldots, \theta^r,
X_1,\ldots,X_s)\nonumber\\
&+\sum_{j=1}^s T(\theta^1,\ldots, \theta^r,
X_1,\ldots,{\mathcal D}X_j,\ldots,X_s),
\end{align}
where we have interpreted $T$ as in \eqref{tensormultilinear}.
This can be proved in an analogous way to the classical case  (see for example \cite[Prop. 2.12 and 2.13]{Oneill}). 

\begin{proposition}
A tensor derivation is determined by its value on $\mathfrak{T}^0_0(M,A)={\mathcal F}(A)$ and on $\mathfrak{X}(M)\subset \mathfrak{T}^1_0(M,A)$.
\end{proposition}
\begin{proof}
As a consequence of the product rule, the tensor derivation  is determined by its value on $\mathfrak{T}^1_0(M,A)$, $\mathfrak{T}^0_1(M,A)$ and $\mathfrak{T}^0_0(M,A)$. Moreover, locally, in a system of coordinates $(\Omega,\varphi)$, any $X\in \mathfrak{T}^1_0(\Omega,A)$ can be expressed as $X=a^i \partial_i$, with $a_i\in {\mathcal F}(T\Omega\cap A)$, which implies, using the product rule, that $({\mathcal D}X)|_{T\Omega\cap A}={\mathcal D}(a^i) \partial_i+a^i {\mathcal D}(\partial_i)$. Finally for $\theta\in \mathfrak{T}^0_1(M,A)$, applying the product rule again, we deduce that
\begin{equation}\label{deroneform}
{\mathcal D}\theta (X)={\mathcal D}(\theta X)-\theta({\mathcal D}(X)).
\end{equation}
\end{proof}
\begin{definition}\label{ander}
 Let us define an {\em anisotropic derivation} $\delta$ in $A$ with associated vector field $Z\in  {\mathfrak X}(M)$ as a map 
 \[\delta: A\times \mathfrak{X}(M)\rightarrow TM,\quad\quad (v,X)\mapsto\delta^v X\in T_{\pi(v)}M,\]
 where, for every $X\in \mathfrak{X}(M)$,  $v\rightarrow \delta^vX$ provides a (smooth) element of $\mathfrak{T}^1_0(M,A)$, and such that
 \begin{enumerate}[(i)]
 \item $\delta^v(X+Y)=\delta^v X+\delta^v Y$, $X,Y\in  {\mathfrak X}(M)$,
\item $\delta^v(fX)=Z(f) X+f \delta^vX $ for any $f\in {\mathcal F}(M)$, $X\in  {\mathfrak X}(M)$,
\end{enumerate}
Moreover, in this case, $(\delta^{V} X)(p):=\delta^{V(p)}X$ is smooth for any $A$-admissble $V\in {\mathfrak X}(\Omega)$, defined in an arbitrary open subset $\Omega\subset M$  and any $X\in  {\mathfrak X}(M)$.
 \end{definition}
 
 \begin{lemma}\label{lemma:dhv}
 Given an anisotropic derivation $\delta$ in $A$ with associated vector field $Z$, we can define the derivation of a function $h\in {\mathcal F}(A)$ as
 \begin{equation}\label{dhv}
 {\mathcal D}(h)(v):=Z(h(V))(\pi(v))-(\partial^\nu h)_v(\delta^vV)
 \end{equation}
 with $V\in{\mathfrak X}(\Omega)$ any $A$-admissible vector field on an open subset $\Omega\subset M$ such that $V(\pi(v))=v$, which satisfies the Leibnitzian property ${\mathcal D}(fg)={\mathcal D}(f)g+f{\mathcal D}(g)$ for any $f,g\in {\mathcal F}(A)$.
 \end{lemma}
 \begin{proof}
The only non-trivial point is to check that ${\mathcal D}(h)(v)$ is well-defined, namely, it does not depend on the extension $V$ of $v$.  Consider a system of coordinates $(\Omega,\varphi)$ and define the non-linear coefficients of $\delta$ as the functions $\delta_{\, \,j}^k(v)$ satisfying that 
 $\delta^v\partial_i= \delta_{\,\, i}^k(v) \partial_k$. 
In the following, we will omit, with an abuse of notation, the evaluation at $v$ and $\pi(v)$. Observe that if, in the system of coordinates, $V=V^k\partial_k$ and $Z= Z^k\partial_k$, then $\delta^vV=Z(V^k) \partial_k+V^j\delta^v(\partial_j)=(Z^j\frac{\partial V^k}{\partial x^j}+V^j\delta_{\,\, j}^k)\partial_k$ and for any $Y\in \mathfrak{T}^1_0(M,A)$, we have $\partial^\nu h(Y)=\frac{\partial h(v+tY)}{\partial t}= Y^k\frac{\partial h}{\partial y^k}.$ Putting this information in \eqref{dhv}, and using the natural coordinates $(T\Omega,\tilde\varphi)$ introduced in \eqref{naturalcoord}, we finally get
 \begin{align*}
 {\mathcal D}(h)&=Z^i\frac{\partial h}{\partial x^i}+Z^i\frac{\partial V^j}{\partial x^i} \frac{\partial h}{\partial y^j}-Z^i\frac{\partial V^k}{\partial x^i} \frac{\partial h}{\partial y^k}-V^j\delta^k_{\,\, j} \frac{\partial h}{\partial y^k}\\
 &=Z^i\frac{\partial h}{\partial x^i}-V^j\delta^k_{\,\, j} \frac{\partial h}{\partial y^k},
 \end{align*}
which concludes that ${\mathcal D}(h)$ does not depend on the considered extension $V$.
 \end{proof}
\begin{remark}
	Observe that using Lemma \ref{lemma:dhv}, one can extend an anisotropic derivation to a derivation acting on $\mathfrak{T}^1_0(M,A)$. By the local nature of derivations, it is enough to define them on $\mathfrak{T}^1_0(\Omega,A)$, with $\Omega$ and open subset of $M$. Indeed, given $X\in \mathfrak{T}^1_0(M,A)$, and a frame $X_1,\ldots, X_n$ of vector fields on $\Omega\subset M$, one has $X=a^iX_i$, with $a^i\in {\mathcal F}(T\Omega\cap A)$ and then
	\begin{equation}\label{defDX} {\mathcal D}(X)(v)={\mathcal D}(a^i)(v)X_i(\pi(v))+a^i(v) \delta^vX_i.
	\end{equation}
	Moreover, it is not difficult to see that it is well-defined, namely, it does not depend on the frame. Indeed, if $Y_1,\ldots,Y_n$ is another frame and $X=b^jY_j$ with $b^j\in {\mathcal F}(T\Omega\cap A)$, we know that $Y_j=c_j^i X_i$, with $c_j^i\in{\mathcal F}(\Omega)$, and  then $a^i=b^jc_j^i$. Observe that we can interpret ${\mathcal F}(\Omega)$ as a subset of ${\mathcal F}(T\Omega\cap A)$ by identifying $f\in {\mathcal F}(\Omega)$ with $f\circ \pi|_{A\cap T\Omega}$. In such a case, for $f\in {\mathcal F}(\Omega)$, with the definition of the above lemma, ${\mathcal D}(f)=Z(f)$, because  $\partial^\nu f=0$. This implies that
		\begin{multline}\label{DXdef} {\mathcal D}(X)(v)={\mathcal D}(a^i)(v)X_i(\pi(v))+a^i(v) \delta^vX_i=
		{\mathcal D}(b^jc_j^i)(v)X_i(\pi(v))+b^j(v)c_j^i \delta^vX_i\\
		=	{\mathcal D}(b^j)(v)c_j^iX_i(\pi(v))+b^j(v) Z(c_j^i)X_i(\pi(v))+b^j(v)c_j^i \delta^vX_i\\={\mathcal D}(b^j)(v)c_j^i X_i(\pi(v))+b^j(v)( Z(c_j^i)X_i(\pi(v))+c_j^i \delta^vX_i)\\={\mathcal D}(b^j)(v)Y_j(\pi(v))+b^j(v)\delta^vY_j,
		\end{multline}
		where we have also used part $(ii)$ in Definition \ref{ander} in the last identity. Moreover, it is straightfoward to check that with this definition
		\[{\mathcal D}(fX)={\mathcal D}(f)X+f{\mathcal D}(X)\]
		for any $f\in {\mathcal F}(A)$. Finally, observe that given $Y\in \mathfrak{T}^1_0(M,A)$ and an $A$-admissible vector field $V\in\mathfrak{X}(\Omega)$ (for a certain neighbourhood $\Omega\subset M$ of $\pi(v)$) that  extends $v\in A$, it follows that 
		\begin{equation}\label{DY}
		{\mathcal D}(Y)(v)=	{\mathcal D}(Y(V))|_v-(\partial^\nu Y)_v(\delta^vV),
		\end{equation}
		where $(\partial^\nu Y)_v(z)=\left.\frac{d}{dt} Y(v+tz)\right|_{t=0}$, for any vector $z\in T_{\pi(v)}M$.
		This can be easily checked using a frame of vector fields as above. In particular, considering a system of coordinates $(\Omega,\varphi)$ and assuming that $Y=a^i\partial_i$, one has
		\begin{align*}
			{\mathcal D}(Y)(v)&=	{\mathcal D}(a^i\partial_i)(v)={\mathcal D}(a^i)|_v\partial_i|_{\pi(v)}+a^i(v) {\mathcal D}(\partial_i)|_v \\
			&=Z(a^i(V))\partial_i|_{\pi(v)}-(\partial^\nu a^i)_v(\delta^vV)\partial_i|_{\pi(v)}+a^i(v) {\mathcal D}(\partial_i)|_v\\
			&={\mathcal D}(a^i(V)\partial_i)|_v-\partial^\nu(a^i\partial_i)_v(\delta^vV)={\mathcal D}(Y(V))|_v-(\partial^\nu Y)_v(\delta^vV).
			\end{align*}
	\end{remark}
 \begin{theorem}
 \label{existderiv}
 Let $\delta$ be an anisotropic derivation in $A\subset TM$ with associated vector field $Z\in {\mathfrak X}(M)$. Then
  there exists a unique tensor derivation $\mathcal D$ such that ${\mathcal D}(X)(v)=\delta^vX$ for every $X\in {\mathfrak X}(M)$, and ${\mathcal D}(h)$ is defined by \eqref{dhv}.
 \end{theorem}
 \begin{proof}

Recall that the derivation of $Y\in \mathfrak{T}^1_0(M,A)$ is well-defined by \eqref{defDX}, and the derivation of a one-form is given by \eqref{deroneform}. Moreover, 
  using \eqref{productrule}, it follows that  the tensor derivation of a tensor $T\in \mathfrak{T}^r_s(M,A)$ is given  by
 \begin{align*}
({\mathcal D}T)(\theta^1,\ldots,\theta^r,X_1,\ldots, X_s)=& 
{\mathcal D}(T(\theta^1,\ldots,\theta^r,X_1,\ldots,X_s))\\
&-\sum_{i=1}^r T(\theta^1,\ldots,{\mathcal D}\theta^i,\ldots, \theta^r,
X_1,\ldots,X_s)\\
&-\sum_{j=1}^s T(\theta^1,\ldots, \theta^r,
X_1,\ldots,{\mathcal D}X_j,\ldots,X_s).
\end{align*}
It is straightforward to check that ${\mathcal D}T$ is ${\mathcal F}(A)$-multilinear and then an anisotropic tensor. 
 Moreover, as in the classical case of tensor derivations (see \cite[Theorem 2.15]{Oneill}), 
 one can check that ${\mathcal D}$ is an anisotropic tensor derivation (it commutes with contractions and it satisfies the Leibnitzian rule for the tensor product).
 \end{proof}
 Let us see how to compute the tensor derivation when we consider the tensor $T$ evaluated in one-forms $\theta^1,\ldots,\theta^r$ and vector fields $X_1,\ldots, X_s$ in $M$ (recall that, by Remark \ref{extendT}, this is enough to determine $T$). In such a case, from the definitions
\begin{multline}\label{formderi}
{\mathcal D}(T(\theta^1,\ldots,\theta^r,X_1,\ldots,X_s))(v)\\=
Z(T_V(\theta^1,\ldots,\theta^r,X_1,\ldots,X_s))(\pi(v))-(\partial^\nu T)_v(\theta^1,\ldots,\theta^r,X_1,\ldots,X_s,\delta^VV),
\end{multline}
where $V$ is an $A$-admissible vector field such that $V(\pi(v))=v$ (recall that this quantity does not depend on the chosen $V$ by Lemma \ref{lemma:dhv}) defined in a neighbourhood of $\pi(v)$.
Observe that in the last expression we are using that 
\[(\partial^\nu T)_v(\theta^1,\ldots,\theta^r,X_1,\ldots,X_s)=(\partial^\nu h)_v(\delta^VV),\] where $h(v)=T_v(\theta^1,\ldots,\theta^r,X_1,\ldots,X_s,\delta^VV)$. Moreover,  because of the tensorial properties, we only need to define the vector fields and the one-forms in a neighbourhood of $\pi(v)$. So from now on we will make an abuse of notation putting together in the evaluation of the anisotropic tensor, one-forms and vector fields defined in different open subsets of $M$ (we can always take the smallest one).
Finally,
\begin{align}\label{completeDer}
({\mathcal D}T)_v(\theta^1,\ldots,\theta^r,X_1,\ldots, X_s):=& 
Z(T_V(\theta^1,\ldots,\theta^r,X_1,\ldots,X_s))(\pi(v))\nonumber\\&-(\partial^\nu T)_v(\theta^1,\ldots,\theta^r,X_1,\ldots,X_s,\delta^VV)\nonumber\\
&-\sum_{i=1}^r T_v(\theta^1,\ldots,{\mathcal D}\theta^i,\ldots, \theta^r,
X_1,\ldots,X_s)\nonumber\\
&-\sum_{j=1}^s T_v(\theta^1,\ldots, \theta^r,
X_1,\ldots,\delta^VX_j,\ldots,X_s).
\end{align}
Moreover, if $X$ is another vector field, ${\mathcal D}\theta^i(X)=Z(\theta^i(X))-\theta^i(\delta^V(X))$.
 \begin{remark}\label{Tensorequiv}
 	Observe that if we consider a tensor as in \eqref{forthecurv}, we can define the derivation of such a tensor by considering the related tensor in \eqref{interpr}, computing its derivation and then going back to the form \eqref{forthecurv} of the derivation. After some straightforward computations, it turns out that the derivation satisfies a formula analogous to \eqref{completeDer}, but replacing the first term $Z(T_V(\theta^1,\ldots,\theta^r,X_1,\ldots,X_s))$ with $\delta^v(T_V(X_1,\ldots,X_s))$. Indeed, if $\bar T(\theta,X_1,\ldots,X_s)=\theta(T(X_1,\ldots,X_s)$, then 
 	\begin{multline*}
 	{\mathcal D}\bar T(\theta,X_1,\ldots,X_s)={\mathcal D}( \bar T(\theta,X_1,\ldots,X_s))- \bar T({\mathcal D}\theta,
 	X_1,\ldots,,X_s)\\-\sum_{j=1}^s \bar T(	\theta,X_1,\ldots,{\mathcal D}X_j,\ldots,X_s),
 	\end{multline*}
 	where $X_1,\ldots,X_s\in  \mathfrak{T}^1_0(M,A)$ and $\theta\in  \mathfrak{T}^0_1(M,A)$.
 	Using \eqref{deroneform} and the definition of $\bar T$, it follows that 
 	\[{\mathcal D}( \bar T(\theta,X_1,\ldots,X_s))- \bar T({\mathcal D}\theta,
 	X_1,\ldots,,X_s)=\theta({\mathcal D}(T(X_1,\ldots,X_s))).
 	\]
 	Putting the last two identities together, we get
 	\begin{equation}\label{derbarT}	{\mathcal D}\bar T(\theta,X_1,\ldots,X_s)=\theta({\mathcal D}(T(X_1,\ldots,X_s)))-\sum_{j=1}^s \bar T(	\theta,X_1,\ldots,{\mathcal D}X_j,\ldots,X_s).
 	\end{equation}
 	We will define
 		\begin{equation}\label{defDT}
 		{\mathcal D}T(X_1,\ldots,X_s)={\mathcal D}(T(X_1,\ldots,X_s))-\sum_{j=1}^s T(X_1,\ldots,{\mathcal D}X_j,\ldots,X_s),
 		\end{equation}
 		 and then, from \eqref{derbarT}, it follows that ${{\mathcal D}\bar T}=\overline{{\mathcal D}T}$, or explicitly, ${{\mathcal D}\bar T}(\theta,X_1,\ldots,X_s)=\theta({\mathcal D}T(X_1,\ldots,X_s)).$  Finally, observe that we can compute ${\mathcal D}(T(X_1,\ldots,X_s)))$, with $X_1\ldots,X_s\in{\mathfrak X}(M)$, using an $A$-admissible vector field $V$ which extends $v\in A$ as in \eqref{DY}, obtaining
 		 \begin{equation}\label{firstDerV}
 		 {\mathcal D}(T(X_1,\ldots,X_s))=\delta^v(T_V(X_1,\ldots,X_s))-
 		 (\partial^\nu T)_v(X_1,\ldots,X_s,\delta^V V).
 		 \end{equation}
 		and extending the vertical derivation in Definition \ref{defi:vert} to a tensor as in \eqref{forthecurv} in the natural way. In order to check this, use a coordinate system $(\Omega, \varphi)$ and express  $X_i=X^{k_i}_i\partial_{k_i}$. Then, for arbitrary vector fields $X_1,\ldots,X_s$,
 		\begin{align*}
 	{\mathcal D}(T(X_1,\ldots,X_s))(v)=&Z(X^{k_1}_1\ldots X^{k_s}_s) T^l_{\,\,k_1\ldots k_s}(v)\partial_l\\&+ X^{k_1}_1\ldots X^{k_s}_s \big(Z(T^l_{\,\,k_1\ldots k_s}(V))\partial_l
 		- (\partial^\nu T^l_{\,\,k_1\ldots k_s})_v(\delta^VV)\partial_l\\&+T^l_{\,\,k_1\ldots k_s}(v)\delta^v\partial_l\\=&\delta^v(T_V(X_1,\ldots,X_s))- (\partial^\nu T^l_{\,\,k_1\ldots k_s})_v(\delta^VV)\partial_l ,
 		\end{align*}
 		where $T^l_{\,\,k_1\ldots k_s}(v)\partial_l=T_v(\partial_{k_1},\ldots,\partial_{k_s})$.
 		It is not difficult to check that 
 		\[(\partial^\nu T^l_{\,\,k_1\ldots k_s})_v(Y)\partial_l=(\partial^\nu T)_v(\partial_{k_1},\ldots,\partial_{k_s},Y),\]
 		for any vector field $Y\in {\mathfrak X}(\Omega)$, and then that \eqref{firstDerV} holds.
 	\end{remark}
\section{Anisotropic linear connections} \label{anlinconn}
\begin{definition}\label{aniconnection}
An  {\em anisotropic (linear) connection} is  a map
 \[\nabla: A\times \mathfrak{X}(M)\times\mathfrak{X}(M)\rightarrow TM,\quad\quad (v,X,Y)\mapsto\nabla^v_XY:=\nabla(v, X,Y)\in T_{\pi(v)}M\]
 such that $\nabla(v,X,\cdot)$ is an anisotropic derivation (recall Definition \ref{ander}) for every $X\in \mathfrak{X}(M)$ with $X$ as associated vector field, and this map is ${\mathcal F}(M)$-linear in $X$, namely,
 \begin{enumerate}[(i)]
 	\item $\nabla^v_X(Y+Z)=\nabla^v_XY+\nabla^v_XZ$, for any $X,Y,Z\in  {\mathfrak X}(M)$,
 	\item $\nabla^v_X(fY)=X(f) Y|_{\pi(v)}+f(\pi(v)) \nabla^v_XY $ for any $f\in {\mathcal F}(M)$, $X,Y\in  {\mathfrak X}(M)$,
 	\item $A\ni v\rightarrow \nabla^v_XY$ belongs to ${\mathfrak T}^1_0(M,A)$ for any $X,Y\in  {\mathfrak X}(M)$,
 	\item $\nabla^v_{fX+hY}Z=f(\pi(v))\nabla^v_XZ+h(\pi(v)) \nabla^v_YZ$, for any $f,h\in {\mathcal F}(M)$, $X,Y,Z\in  {\mathfrak X}(M)$.
 	\end{enumerate}
 \end{definition}
  We define the {\em torsion} of $\nabla$ as
\begin{equation}\label{torsion}
{\mathcal T}_v(X,Y)=\nabla^v_XY-\nabla^v_YX-[X,Y].
\end{equation}
 It is not difficult to prove that $\mathcal T$ is an anisotropic tensor (recall Remark \ref{extendT}). One should have in mind that it is not possible to evaluate \eqref{torsion} in elements of $\mathfrak{T}^1_0(M,A)$ rather than those of ${\mathfrak X}(M)$. This is because in such a case $[X,Y]$ is not well-defined. But it can be extended to $\mathfrak{T}^1_0(M,A)\times \mathfrak{T}^1_0(M,A)$ as in Remark \ref{extendT}.

  Moreover, we will say that the anisotropic linear connection is {\em  torsion-free} if ${\mathcal T}=0$ or, equivalently,
 \[\nabla^v_XY-\nabla^v_YX=[X,Y]\]
 for every $X,Y\in \mathfrak{X}(M)$ and $v\in A$. 
Given a system of coordinates $(\Omega,\varphi)$, we define the {\em Christoffel symbols} of $\nabla$ as the functions $\Gamma_{\,\,jk}^i:A\cap T\Omega\rightarrow \R$ satisfying that $\nabla^v_{\partial_j}\partial_k=\Gamma_{\,\, jk}^i(v) \partial_i$.
 \begin{remark}\label{torsionfree}
  Observe that an anisotropic connection is torsion-free if and only if its Christoffel symbols are symmetric in $j$ and $k$.
  \end{remark}
\begin{remark}\label{derT}Applying  Theorem \ref{existderiv}, we can get the derivative of any tensor $T\in \mathfrak{T}^r_s(M,A)$ with respect to a vector field $X\in \mathfrak{X}(M)$ obtaining a new tensor $\nabla_XT  \in \mathfrak{T}^r_s(M,A)$. Moreover, we can also define an anisotropic tensor $\nabla T  \in \mathfrak{T}^{r}_{s+1}(M,A)$, with $\nabla T(\dots,X)=\nabla_XT(\dots)$. As usual, we can extend the evaluation in the last component to $\mathfrak{T}^1_0(M,A)$ as in Remark \ref{extendT}.  Regarding Remark \ref{Tensorequiv}, when a tensor is interpreted as in \eqref{forthecurv}, we can also define a derivation $\nabla_X$ as in \eqref{defDT}. In such a case, if we fix an $A$-admissible vector field $V$, then the second term on the right hand side in \eqref{defDT} can be computed as
\begin{multline}\label{partialTexpres}
{\mathcal D}(T(X_1,\ldots,X_s))(v)=\nabla_X(T(X_1,\ldots,X_s))(v)\\=\nabla^V_X(T_V(X_1,\ldots,X_s))(\pi(v))-(\partial^\nu T)_v(X_1,\ldots,X_s,\nabla^V_XV),
\end{multline}
where $X_1,\ldots, X_s\in {\mathfrak X}(M)$ and $\partial^\nu T$ is a natural generalization of Definition \ref{partialT} to the tensors as in \eqref{forthecurv} (see \eqref{firstDerV}). 
\end{remark}
 
 \subsection{Associated non-linear  and Ehresmann connections}\label{nonlinearity}  Observe that an anisotropic connection determines a {\it non-linear connection} as follows:
 \[D:  \mathfrak{X}(M)\times\mathfrak{X}^A(\Omega)\rightarrow TM,\quad\quad (X,Y)\mapsto D_XY:=\nabla^Y_XY,\]
 where $\mathfrak{X}^A(\Omega)$ is the subset of $A$-admissible vector fields in an open subset $\Omega\subset M$. It is not difficult to check that
 $D$ is linear in $X$ but not in $Y$. 
 On the other hand, a couple of vector fields $(X,Y)\in \mathfrak{X}(M)\times\mathfrak{X}^A(\Omega)$ determine a vector $W$ in $T_{Y(p)}(TM)$, namely, if $\alpha$ is the  integral curve of $X$ through $\pi(Y(p))$, $W$ is the velocity of the restriction of $Y$ to $\alpha$, or in natural  coordinates $(T \Omega,\tilde{\varphi})$ for $TM$ with $\tilde{\varphi}=(x^1,\ldots,x^n,y^1,\ldots,y^n)$ associated with a coordinate system $(\Omega, \varphi=(x^1,\ldots,x^n))$ of $M$, 
 \[W(Y)=(X^1,\ldots, X^n,X(Y^1),\ldots,X(Y^n)).\] Moreover, $D_XY$ only depends on $W$ at each point, since in coordinates, we have that
 \[D_XY=(X(Y^k)+X^i Y^j \Gamma_{\, \,ij}^k(Y))\partial_k.\] We will denote by $N_{\,\, i}^k(v)=v^j \Gamma_{\,\, ij}^k(v)$, with $v=v^i\partial_i|_{\pi(v)}$, the {\em non-linear coefficients} associated with $\nabla$. 
 \begin{definition}
 The {\em Ehresmann connection} associated with the non-linear connection is given by the horizontal space obtained by the vector fields $W$ in $TTM$ such that $D_XY=0$, or, equivalently,  if $W=W^i_1\frac{\partial}{\partial x^i}+W^j_2 \frac{\partial}{\partial y^j}$ in a natural coordinate system $(T\Omega,\tilde\varphi)$, then $W\in T_v(TM)$ with $v\in A$  belongs to the horizontal space if and only if
 \[W_2^k+W_1^i N_{\,\, i}^k(v)=0,\quad k=1,\ldots,n.\]
 \end{definition}
  A basis of the horizontal space in $v\in A$ is given by $\frac{\delta}{\delta x^i}|_v=\frac{\partial}{\partial x^i}|_v-N^k_{\,\, i}(v) \frac{\partial}{\partial y^k}|_v, i=1,\ldots, n$. 
  
  Let us obtain coordinates formulae for the tensor derivatives using the non-linear connection. As a first observation, when se consider ${\mathcal D}=\nabla_{\partial_i}$, $i\in\{1,\ldots,n\}$, in Lemma \ref{lemma:dhv} and $h\in {\mathcal F}(A)$, then
  \begin{equation}\label{nonlinearh}
  \nabla_{\partial_i} (h)=\frac{\partial h}{\partial x^i}-N^k_{\,\, i} \frac{\partial h}{\partial y^k}=\frac{\delta}{\delta x^i}(h):=\frac{\delta h}{\delta x^i}. 
  \end{equation}
  This follows easily from the last formula in the proof of Lemma \ref{lemma:dhv}, since in this case, $Z=\partial_i$ and $\delta^k_{\,\,j}=\Gamma^k_{\,\, ij}$.
   Observe that using \eqref{deroneform}, one can show that $\nabla^v_{\partial_i} (dx^j)=-\Gamma^j_{\,\, ik}(v) dx^k|_{\pi(v)}$, and by direct computations, using part $(ii)$ in the definition of anisotropic connection and \eqref{nonlinearh},
 we deduce that 
 \begin{equation}
 T^i_{\,\,j|k}(v):=(\nabla_{\partial_k}T)_v(dx^i,\partial_j)=\frac{\delta T^i_{\,\, j}}{\delta x^k}(v)+\Gamma^i_{\,\, kl}(v) T^l_{\,\,j}(v)-\Gamma^l_{\,\, kj}(v) T^i_{\,\,l}(v),
 \end{equation} 
 where $\frac{\delta T^i_{\,\, j}}{\delta x^k}(v)=\frac{\partial T^i_{\,\, j}}{\partial x^k}(v)-N^l_{\,\, k}(v)\frac{\partial T^i_{\,\, j}}{\partial y^l}(v)$.
 Observe that the expression for $T^i_{\,\, j|k}$ is what in many references call ``the horizontal derivative''.
\subsection{Curvature tensors of the anisotropic connection} \label{curvanis}
Given an anisotropic (linear) connection $\nabla$, we can define the associated curvature tensor $R_v:{\mathfrak X}(M)\times {\mathfrak X}(M)\times {\mathfrak X}(M)\rightarrow T_{\pi(v)}M$, as follows
\begin{equation}\label{Rv}
R_v(X,Y)Z=\nabla^v_X(\nabla_YZ)-\nabla^v_Y(\nabla_XZ)-
\nabla^v_{[X,Y]}Z,
\end{equation}
for any $v\in A$ and $X,Y,Z\in {\mathfrak X}(M)$.
In this expression, let us point out that we interpret $\nabla_YZ$ as an element of $\mathfrak{T}^1_0(M,A)$, in the sense that $\nabla_YZ(v)=\nabla^v_YZ$, and $\nabla^v_X$ as the extension of ${\mathcal D}=\nabla^v_X$ described in \eqref{DXdef}. It is not difficult to prove that $R_v$ is ${\mathcal F}(M)$-multilinear, and then it determines an anisotropic tensor (recall Remark \ref{extendT}). As in the case of the torsion, $R_v$ cannot evaluated directly in elements of $\mathfrak{T}^0_1(M,A)$, as the Lie bracket is not well-defined, but it can be extended by linearity.

Let us give the expression of the curvature in coordinates. If
$R_v(\partial_k,\partial_l)\partial_j=R_{j\,\, kl}^{\,\,\,i}(v)\partial_i$, then from a direct computation using \eqref{nonlinearh} and the properties of $\nabla^v$, it follows that
\[R_{j\,\, kl}^{\,\,\,i}=\frac{\delta \Gamma^i_{\,\, lj}}{\delta x^k}-\frac{\delta \Gamma^i_{\,\, kj}}{\delta x^l}+\Gamma^i_{\,\, kh}\Gamma^h_{\,\, lj}-
\Gamma^i_{\,\, lh}\Gamma^h_{\,\, kj}.
\]
\subsection{Covariant derivatives along curves}\label{covariant}
In the following, given a smooth curve $\gamma:[a,b]\rightarrow M$, $\mathfrak{X}(\gamma)$ will denote the space of smooth vector fields along $\gamma$ and ${\mathcal F}(I)$ the smooth real functions defined on $I=[a,b]$.
\begin{definition}
	An {\em anisotropic covariant derivation} $D^v_\gamma$ in $A$ along a curve $\gamma:[a,b]\rightarrow M$ is a map 
	\[D^v_\gamma: \mathfrak{X}(\gamma)\rightarrow T_{\pi(v)}M,\quad\quad X\mapsto D_\gamma^vX \]
	for every $v\in A$ with $\pi(v)=\gamma(t_0)$, and $t_0\in [a,b]$, such that 
	\begin{enumerate}[(i)]
		\item $D^v_\gamma(X+Y)=D^v_\gamma X+D^v_\gamma Y$, $X,Y\in  {\mathfrak X}(\gamma)$,
		\item $D_\gamma^v(fX)=\frac{df}{dt}(t_0) X(t_0)+f(t_0) D_\gamma^vX $ for any $f\in {\mathcal F}(I)$, $X\in  {\mathfrak X}(\gamma)$,
		\item $D_\gamma^VX(t):=D_\gamma^{V(t)}X$ is smooth for any $V\in {\mathfrak X}(\gamma)$, $X\in  {\mathfrak X}(\gamma)$ and $V$, $A$-admissible, that is, $V(t)\in A$ for every $t\in [a,b]$.
	\end{enumerate}
\end{definition}
\begin{proposition}
	Every anisotropic linear connection $\nabla$ determines an induced anisotropic covariant derivative along a smooth curve $\gamma:[a,b]\rightarrow M$, with the following property: if $X\in {\mathfrak X}(M)$, then $D_\gamma^v (X_{\gamma})=\nabla^v_{\dot\gamma} X$, where $X_\gamma$ is the vector field in ${\mathfrak X}(\gamma)$ defined as $X_\gamma(t)=X(\gamma(t))$ for every $t\in [a,b]$.  
\end{proposition}
\begin{proof}
	It follows as in \cite[Proposition 3.18]{Oneill}. Observe that  in coordinates around $\gamma(t_0)=\pi(v)$, with $t_0\in [a,b]$, the induced covariant connection is given by
	\begin{equation}
	D_\gamma^v X=\dot X^i(t_0)\frac{\partial }{\partial x^i}+X^i(t_0)\dot\gamma^j(t_0) \Gamma^k_{\, ij}(v) \frac{\partial}{\partial x^k}.
	\end{equation}
\end{proof}
\begin{definition}
	We say that a smooth curve $\gamma:[a,b]\rightarrow M$ is $A$-admissible if $\dot\gamma\in A$ for all $t\in [a,b]$. Moreover, we say that an $A$-admissible smooth curve is an {\em autoparallel} of the anisotropic linear connection $\nabla$ if $D^{\dot\gamma}_\gamma \dot\gamma=0$, where $D_\gamma$ is the anisotropic covariant derivative associated with $\nabla$.
\end{definition}
In coordinates, autoparallel curves are given by the equation
\begin{equation}
\ddot \gamma^k+\dot\gamma^i\dot\gamma^j \Gamma^{k}_{\, ij}(\dot\gamma)=0.
\end{equation}
\section{Classical connections as anisotropic connections}

\subsection{Chern connection as Levi-Civita connection}\label{chernconnection}
We say that a smooth, positive two-homogeneous function $L:A\subset TM\setminus 0\rightarrow \R$ is a pseudo-Finsler metric if $A$ is a conic subset, namely, for every $v\in A$ and $\lambda>0$ we have that $\lambda v\in A$ and, in addition, the fundamental tensor defined as 
\begin{equation}
g_v(u,w):=\frac 12 \frac{\partial^2}{\partial t\partial s}L(v+tu+sw)|_{t=s=0}
\end{equation}
for every $v\in A$ and $u,w\in T_{\pi(v)}M$, is non-degenerate. Moreover, we define the Cartan tensor associated with $L$ as
\begin{equation}
C_v(w_1,w_2,w_3):=\left.\frac 14 \frac{\partial^3}{\partial s_3\partial s_2\partial s_1}L\left(v+\sum_{i=1}^3s_iw_i\right)\right|_{
s_1=s_2=s_3=0}.
\end{equation}
Recall that $C_v$ is symmetric and, by homogeneity, one has that $C_v(v,u,w)=C_v(u,v,w)=C_v(u,w,v)=0$ for any $v\in A$ and $u,w\in T_{\pi(v)}M$.
Observe that given any anisotropic linear connection we can compute the covariant derivative of a pseudo-Finsler metric as
\[(\nabla_XL)_v=X(L(V))(\pi(v))-(\partial^\nu L)_v(\nabla^V_XV),\]
where $V$ is any local $A$-admissible extension of $v$ and $X$ is any vector field. Taking into account that the vertical derivative of $L$ satisfies $(\partial^\nu L)_v(w)=2g_v(v,w)$, 
we conclude that 
\[(\nabla_XL)_v=X(L(V))(\pi(v))-2g_v(V,\nabla^V_XV).\]
The covariant derivative of its fundamental tensor can be computed as
\[(\nabla_X g)_v(Y,Z)=X(g_V(Y,Z))(\pi(v))-g_v(\nabla^V_XY,Z)-g_v(Y,\nabla^V_XZ)-(\partial^\nu g)_v(Y,Z,\nabla^V_XV),\]
for $X,Y,Z\in {\mathfrak X}(M)$ and $V$, an $A$-admissible local extension of $v$.
As the vertical derivative of $g$ satisfies $(\partial^\nu g)_v(X,Y,Z)=2C_v(X,Y,Z)$, then
\[(\nabla_X g)_v(Y,Z)=X(g_V(Y,Z))(\pi(v))-g_v(\nabla^V_XY,Z)-g_v(Y,\nabla^V_XZ)-2C_v(Y,Z,\nabla^V_XV).\]
This implies straightforwardly that the Chern connection (when interpreted as a family of affine connections) is the only anisotropic linear connection which is torsion-free and compatible with the metric in the sense that $\nabla g=0$, namely, the Chern connection is the Levi-Civita connection of the pseudo-Finsler metric (see \cite{Mat80}, \cite[Eqs. (7.20) and (7.21)]{Sh01} or \cite[Proposition 2.3]{Jav14a}).  It is possible to get a Koszul formula for this connection:
\begin{multline*}
 2 g_v(\nabla^V_XY,Z)= (X (g_V(Y,Z))-Z (g_V(X,Y))+Y (g_V(Z,X)))(\pi(v))\\
+g_v([X,Y],Z)+g_v([Z,X],Y)-g_v([Y,Z],X)\\
2( -C_v(Y,Z,\nabla^V_XV)-C_v(Z,X,\nabla^V_YV)+C_v(X,Y,\nabla^V_ZV)).
\end{multline*}
\subsection{Nonlinear connection of a spray}\label{nonlinear}
 Let us introduce some general notation. Given a manifold $M$, we will denote by $\pi_M:TM\rightarrow M$ the natural projection from its tangent bundle. In particular, $\pi_{TM}:TTM\rightarrow TM$ denotes the natural projection from the tangent bundle $TTM$ of $TM$. Recall that the vertical distribution of $TM$ is defined as ${\mathcal V}TM=\ker (\pi_{TM})$ (we will also use ${\mathcal V}_vTM={\mathcal V}TM\cap \pi_{TM}^{-1}(v)$). Moreover, we will denote by $\C$ the {\em canonical vertical vector field } defined as $\C(v)=i_v(v)$, where
 \begin{equation}\label{vertlift}
 i_v: T_{\pi(v)}M\rightarrow {\mathcal V}_v TM
 \end{equation}
 is given by $i_v(u)=\frac{d}{dt}(v+ tu)|_{t=0}$.
 \begin{definition}
 Given a manifold $M$, a {\em semi-spray } on $M$ is a vector field  $S$  in a  subset $A\subset TM$   which has the following property: 
\begin{equation}\label{semispray} 
\text{ if $\beta: [a,b]\subset \R\rightarrow A$ is an integral curve of $S$, then $\beta=\dot\alpha$, where $\alpha=\pi_M\circ \beta$. }
\end{equation}
\end{definition} 
 Condition \eqref{semispray} is equivalent to $\pi_{TM}(S)=d\pi_M(S)$. Moreover,  
 in natural coordinates  $(T\Omega,\tilde{\varphi})$ for $TM$ with $\tilde{\varphi}=(x^1,\ldots,x^n,y^1,\ldots,y^n)$ associated with a coordinate system $(\Omega,\varphi=(x^1,\ldots,x^n))$ of $M$, $S$ is a semi-spray if and only if it is expressed as
\begin{equation}\label{Sspray}
S(v)=y^i \frac{\partial}{\partial x^i}- 2G^i(v)\frac{\partial}{\partial y^i},
\end{equation}
where $G^i:T\Omega\cap A\rightarrow \R$ are the so-called {\it coefficients of the semi-spray} and $v=\tilde\varphi^{-1}(x,y)$. 
\begin{definition}
We say that a semi-spray $S$ is a {\it spray} if, in addition, $A$ is conic, namely, if $v\in A$, then $\lambda v\in A$ for every $\lambda>0$ and $S$ has the following property: 
\begin{multline}\label{spray}
\text{If $\beta=\dot\alpha$ is an integral curve of $S$, then $\tilde{\beta}(t)=\lambda \dot\alpha(\lambda t)$}\\
\text{ is an integral curve of $S$ for every $\lambda>0$.}
\end{multline}
\end{definition}
Condition \eqref{spray} is equivalent to the  positive homogeneity of degree two of the functions $G^i$ (recall \eqref{Sspray}) on $T\Omega\cap A$ (see \cite[Chapter 4]{Sh01}). It is also equivalent to the equation $[\C,S]=S$, where $\C$ is the canonical vertical vector field. We will say that  the {\it autoparallel curves of the spray} are the projections in $M$ of the integral curves of $S$ in $A$.  In a system of coordinates $(\Omega,\varphi)$, a curve $\gamma:[a,b]\rightarrow \Omega\subset M$ is an autoparallel curve of the spray $S$ if and only if
\begin{equation}
\ddot \gamma^i=- 2G^i(\dot \gamma),
\end{equation}
for $i=1,\ldots,n$, where $\gamma^i$ are the coordinates of $\gamma$.
Moreover, a spray determines naturally an Ehresmann connection. A basis of the horizontal space is given by
\[\frac{\delta}{\delta x^j}=\frac{\partial}{\partial x^j}-N_{\,j}^i(v) \frac{\partial}{\partial y^i},\quad j=1,\ldots,n\]
where $N^i_{\,j}(v)=\frac{\partial G^i}{\partial y^j}(v)$ are called the non-linear coefficients. Having a horizontal subspace, we have a splitting of $TTM$:
\[TTM= {\mathcal V}TM+{\mathcal H} TM,\]
which allows us to define the vertical and the  horizontal projections.  We can also define the non-linear covariant derivative of any vector $X$ along a curve $\alpha:[a,b]\rightarrow M$ as $D_\alpha X=i_v^{-1}({\mathcal V} \dot X)$. In a  system of coordinates,
\[D_\alpha X = \dot X^i \partial_i+ \dot\alpha^i N^j_{\,\, i}(X) \partial_j,\]
where $X^i$ are the coordinates of $X$ in $(\Omega,\varphi)$. Then it follows from the homogeneity of $G^i$ that $\gamma:[a,b]\rightarrow M$ is an autoparallel curve of the spray if and only if $D_\gamma \dot\gamma=0$, namely, the nonlinear covariant derivative of the velocity is zero. Indeed, the homogeneity implies that $2G^j(X)=X^i N^j_{\,\, i}(X)$, and then $\dot\gamma^i N^j_{\,\, i}(\dot\gamma)=2G^j(\dot\gamma)$.

\subsection{Berwald connection and the classical tensors}\label{berwaldcon}
We have seen in subsection \ref{nonlinear} how to associate a non-linear connection with a spray. Moreover, autoparallel curves can be recovered as the curves having parallel velocity with the non-linear connection, namely, in a natural coordinate system $(T\Omega,\tilde{\varphi})$, they satisfy
\[\ddot \gamma^i=-\dot \gamma^k N^i_{\,\, k}(\dot\gamma),\]
since, by homogeneity, we can apply Euler's theorem to obtain that 
\[y^k N^i_{\,\, k}(\tilde\varphi^{-1}(x,y))=2 G^i(\tilde\varphi^{-1}(x,y)),\quad i=1,\ldots, n.\] If we want to find an anisotropic connection with the same autoparallel curves as the spray, then we need that 
\[ \dot\gamma^i N^k_{\,\, i}(\dot\gamma)= \dot \gamma^i \dot \gamma^j \Gamma^k_{\,\, ij}(\dot\gamma),\]
being 
\begin{equation}\label{ChrisBer}
\Gamma^k_{\,\, ij}(v)=\frac{\partial N^k_{\,\,i}}{\partial y^j}(v)=\frac{\partial^2 G^k}{\partial y^i\partial y^j}(v)
\end{equation} 
a possible choice, since applying Euler's theorem for (positive) homogeneous functions to $\frac{\partial G^k}{\partial y^i}(\tilde\varphi^{-1}(x,y))$, which is positive homogeneous of degree one in $y$, it follows that $\frac{\partial G^k}{\partial y^i}(\tilde\varphi^{-1}(x,y))=y^j \frac{\partial^2 G^k}{\partial y^j\partial y^i}(\tilde\varphi^{-1}(x,y))$, and then
\[y^iy^j \frac{\partial^2 G^k}{\partial y^i\partial y^j}(\tilde\varphi^{-1}(x,y))=y^i \frac{\partial G^k}{\partial y^i}(\tilde\varphi^{-1}(x,y))=y^i N^k_{\,\, i}(\tilde\varphi^{-1}(x,y)). \]
 The quantities in \eqref{ChrisBer} are the Christoffel symbols of the Berwald connection (see \cite[\S 7]{Sh01}\ and \cite[\S 8]{HJP15} for a treatment of the Berwald connection as an anisotropic linear connection). In the following, we will denote the Chern  and Berwald connections, respectively, as  $\nabla$ and $\tilde{\nabla}$.

Let us now introduce some classical tensors. 
Recall that a pseudo-Finsler metric $L:A\subset TM\setminus 0\rightarrow \R$ determines a spray with
\begin{equation}\label{sprayfinsler}
G^i(\tilde\varphi^{1}(x,y))=y^jy^k g^{is}\frac{1}{4} \left( 2\frac{\partial g_{sj}}{\partial x^k}-
\frac{\partial g_{jk}}{\partial x^s}\right),
\end{equation}
for a system of natural coordinates $(T\Omega,\tilde{\varphi})$ (see for example \cite[Eq. (4.30)]{Sh01}), where $g_{ij}(v)=g_v(\partial_i,\partial_j)$ and $\{g^{ij}\}$ is the inverse matrix of $\{g_{ij}\}$.
\begin{definition}
Let $S$ be a spray in $A\subset TM\setminus 0$, $\tilde\nabla$ its Berwald connection, $v\in A$, $u,w,z\in T_{\pi(v)}M$ and $V,X,Y$ and $Z$ are arbitrary local extensions of $v,u,w$ and $z$, being $V$, $A$-admissible.
The {\em Berwald tensor} $B$ is defined as the vertical derivative of the Berwald connection, namely,
\begin{equation}\label{berwaldtensor}
B_v(u,w,z)= \frac{\partial}{\partial t}\left(\tilde{\nabla}^{V+tZ}_XY\right)|_{t=0}
\end{equation}
 (see (6.4) in \cite{Sh01}).  
Assume now that $(M,L)$ is a pseudo-Finsler manifold, $S$ is the spray associated with it (see \eqref{sprayfinsler}) and $\nabla$, the Chern connection associated with $L$ (see \S \ref{chernconnection}).
We define the {\it Chern tensor}  of $L$ as the vertical derivative of the Chern connection
\begin{equation}\label{cherntensor}
P_v(u,w,z)= \frac{\partial}{\partial t}\left(\nabla^{V+tZ}_XY\right)|_{t=0}
\end{equation}
 (see (7.23) in \cite{Sh01}, where it has the opposite sign). 
As $\tilde{\nabla}$ and $\nabla$ are torsion-free, $B$ and $P$ are symmetric in the first two components, and by homogeneity, it follows that $B_v(u,w,v)=P_v(u,w,v)=0$. Moreover, the Berwald tensor is symmetric, since in coordinates is given by $B_v(\partial_i,\partial_j,\partial_k)=B_{\,ijk}^l(v)\partial_l$ with $B_{\,ijk}^l=\frac{\partial^3 G^l}{\partial y^i\partial y^j\partial y^k}(v)$ and then
\begin{equation}\label{Berv0}
B_v(v,u,w)=B_v(u,v,w)=B_v(u,w,v)=0.
\end{equation}
Finally, we define the {\em  Landsberg curvature} of a pseudo-Finsler metric $L$ as
 \begin{equation}\label{landsberg}
 \mathfrak{L}_v(u,w,z)=\frac{1}{2}g_v(B_v(u,w,z),v)
 \end{equation}
 (see (6.25) in \cite[Definition 6.2.1]{Sh01} where it has a different sign). From \eqref{Berv0}, it follows that 
 \begin{equation}\label{Land0}
\mathfrak{L}_v(v,u,w)=\mathfrak{L}_v(u,v,w)=\mathfrak{L}_v(u,w,v)=0.
\end{equation}
 \end{definition}
With these definitions, we can write down the difference tensor between the Chern and Berwald connection as 
\begin{equation}\label{differtensor}
 \nabla^v_XY-\tilde{\nabla}^v_XY=\mathfrak{L}^\flat_v(X,Y),
\end{equation}
where  $\mathfrak{L}^\flat$ is determined by $g_v(\mathfrak{L}^\flat_v(u,w),z)=\mathfrak{L}_v(u,w,z)$ (see (7.17) in \cite{Sh01} and observe that the notation for the Chern and Berwald connections is changed). 
\subsection{Relation with linear connections in fiber bundles}\label{relationbundles}

Given a linear connection $\nabla$ in the fiber bundle $\pi_{TM}: {\mathcal V}TM\rightarrow A\subset TM$ and a horizontal connection $\mathcal H$ on $\pi_{TM}:TTM\rightarrow A$, we can define an associated anisotropic connection. First, observe that the horizontal connection is defined as a smooth choice of a $n$-dimensional transverse subspace to ${\mathcal V}_vTM$ at every $v\in A$. This is determined by the lifts to the horizontal subspace of the partial vectors $\partial_i$, $i=1,\ldots,n$, in a coordiante system $(\Omega,\varphi)$, namely,
\[\left.\frac{\delta}{\delta x^i}\right|_v:=\left.\frac{\partial}{\partial x^i}\right|_v-N^j_{\,\,i}(v) \left.\frac{\partial}{\partial y^j}\right|_v.\]
We will assume that the non-linear coefficients are (positive) homogeneous of degree one.
 Recall the map $i_v$ in \eqref{vertlift}. We can define the vertical lift in $v\in A$ of a vector $u\in T_{\pi(v)}M$ as $u^{\mathcal V}=i_v(u)$. Moreover, we can also define the horizontal lift $u^{\mathcal H}$ in $v\in A$ of a vector $u\in T_{\pi(v)}M$ as the unique horizontal vector in $T_vTM$ such that  $\pi_{TM}(u^{\mathcal H})=u$. In coordinates, if $X=a^i \frac{\partial}{\partial x^i}$, then
\[X^{\mathcal H}= a^i \left(\frac{\partial}{\partial x^i}-N^j_{\,\,i} \frac{\partial}{\partial y^j}\right).\]
Observe that using the map $i_v$ we can identify the fiber bundle $\pi_{TM}: {\mathcal V}TM\rightarrow A\subset TM$ with the pullback one introduced in \S \ref{antensor}, $\pi^*_A:\pi^*_A(M)\rightarrow A$. Therefore, a classical connection is a map
\[\nabla: {\mathfrak X}(A)\times {\mathfrak T}_0^1(M,A)\rightarrow {\mathfrak T}_0^1(M,A),\]
which is ${\mathcal F }(A)$-multilinear in the first component, and it is linear and it satisfies the Leibnitz rule in the second one. Using the horizontal connection one has that $X={\mathcal V}(X)+{\mathcal H}(X)$. Moreover, we have the following identifications:
\begin{align*}
i_v:  T_{\pi(v)}M\rightarrow {\mathcal V}(T_v(TM)),\\
\pi_{TM}: {\mathcal H}(T_v(TM))\rightarrow T_{\pi(v)}M.
\end{align*}
Observe that the map $(\cdot)^{\mathcal V}$ is the first map, while $(\cdot)^{\mathcal H}$ is the inverse of the above restriction of $\pi_{TM}$. Moreover, given $X\in{\mathfrak X}(A)$, we can define two elements in ${\mathfrak T}_0^1(M,A)$, namely,
\begin{align*}
A\ni v\rightarrow i_v^{-1}({\mathcal V}(X))\in T_{\pi(v)}M,\\
A\ni v\rightarrow \pi_{TM}({\mathcal H}(X)|_v)\in T_{\pi(v)}M.
\end{align*}
The {\em anisotropic connection associated with the classical connection $\nabla$ } is defined as
\begin{equation}\label{anicon}
\hat\nabla^v_XY= (\nabla_{X^{\mathcal H}} Y)(v),
\end{equation}
for any $v\in A$ and $X,Y\in {\mathfrak X}(M)$. Here we have interpreted $Y$ as an element in $\mathfrak{T}^1_0(M,A)$ as  we consider the natural inclusion ${\mathfrak X}(M)\subset \mathfrak{T}^1_0(M,A)$, namely, we identify $Y$ with $Y^{\mathcal V}$.  On the other hand, observe that given a vector field $X\in {\mathfrak X}(M)$, the horizontal connection also defines a derivative for a function $f\in {\mathcal F}(A)$ as ${\mathcal D}(f)=X^{\mathcal H}(f)$. In principle, this derivative does not have to coincide with the one defined in Lemma \ref{lemma:dhv}, but from now on, we will assume that this is the case. In such a case, in a natural system of coordinates $(T\Omega,\tilde\varphi)$, it holds
\begin{equation}\label{Ngamma}
N^k_{\,\, l}(v)=v^j \Gamma_{\,\, lj}^k(v)\quad \text{for $v\in A\cap T\Omega$}.
\end{equation}
We can extend the anisotropic connection to elements $Z\in  \mathfrak{T}^1_0(M,A)$ as in \eqref{defDX} and, in such a case, it is still true that $\hat\nabla_XZ=\nabla_{X^{\mathcal H}} Z$ for any $v\in A$,  $X\in {\mathfrak X}(M)$ and $Z\in \mathfrak{T}^1_0(M,A)$.
Moreover, we can define an anisotropic tensor
\begin{equation}\label{tensorCb}
\Cb_v(X,Y)=(\nabla_{X^{\mathcal V}} Y)(v),
\end{equation}
for any $v\in A$ and $X,Y\in {\mathfrak X}(M)$.
Observe that $\Cb$ is a tensor because if $f\in {\mathcal F}(M)$, then $Y^{\mathcal V}(f\circ \pi)=0$, but if one considers $Z\in  \mathfrak{T}^1_0(M,A)$, then 
\begin{equation}\label{nablaXVZ}
\nabla_{X^{\mathcal V}} Z=(\partial^\nu Z)_v(X)+\Cb_v(X,Z),
\end{equation}
for any $v\in A$ and  $X\in {\mathfrak X}(M)$. Here we have extended the tensor $\Cb_v$ to $ \mathfrak{T}^1_0(M,A)\times  \mathfrak{T}^1_0(M,A)$, using $\mathcal F(A)$-multilinearity (recall Remark \ref{extendT}).
\begin{proposition}
A classical connection $\nabla$ in the vertical bundle and a horizontal connection on $\pi_{TM}:{\mathcal V}TM\rightarrow A$ are determined by the anisotropic connection $\hat\nabla$ in \eqref{anicon} and the anisotropic tensor $\Cb\in \mathfrak{T}_2^0(M,A)$ in \eqref{tensorCb}, whenever \eqref{Ngamma} holds. The converse also holds, namely, given an anisotropic connection $\hat\nabla$ and an anisotropic tensor $\Cb\in \mathfrak{T}_2^0(M,A)$, we can construct a classical linear connection on the fiber bundle $\pi: {\mathcal V}TM\rightarrow A\subset TM$ associated with the spray determined by the anisotropic connection using the non-linear connection associated with $\hat \nabla$ (see \S \ref{nonlinearity}). 
\end{proposition}
\begin{proof}
	The first claim follows from \eqref{anicon} and \eqref{nablaXVZ}, together with the fact that the derivation that $\hat\nabla_X$ induces on ${\mathcal F} (A)$ is just the one provided by $X^{\mathcal H}$ (see \eqref{nonlinearh}). For the second claim, just define the linear connection
	\[\nabla_XZ=\nabla_{{\mathcal V}(X)}Z+\nabla_{{\mathcal H}(X)}Z=(\partial^\nu Z)_v(i_v^{-1}({\mathcal V}(X)))+\Cb_v(i_v^{-1}({\mathcal V}(X)),Z)+\hat\nabla^v_{\pi_{TM}({\mathcal H}(X))}Z,\]
   being $Z\in \mathfrak{T}^1_0(M,A)$ and $X\in {\mathfrak X}(A)$. It is straightforward to check that the above $\nabla$ satisfies the properties of a classical linear connection.
	\end{proof}
In the following, we will consider a linear connection $\nabla$ on the vertical bundle, a horizontal connection with coefficients $N^j_{\,\, i}(v)$ in a coordinate system, the associated anisotropic connection $\hat \nabla$ satisfying \eqref{Ngamma} and the  anisotropic tensor $\Cb$ given in \eqref{tensorCb}.
\begin{proposition}\label{relationcurder}
Given an arbitrary anisotropic tensor $T$ of the fiber bundle  $\pi_{TM}: {\mathcal V}TM\rightarrow A\subset TM$, we have that $\nabla_{Z^{\mathcal H}}T=\hat\nabla_ZT$ for any vector field $Z\in {\mathfrak X}(M)$ and
\begin{multline}\label{nablaZv}
(\nabla_{Z^{\mathcal V}}T)(\theta^1,\ldots, \theta^r,X_1,\ldots,X_r)(v)=(\partial^\nu T)_v(\theta^1,\dots,\theta^r,X_1,\ldots,X_s,Z) \\
+\sum_{i=1}^rT_v(\theta^1,\ldots,\theta^i(\Cb_v(Z,\cdot)),\ldots,\theta^r,X_1,\ldots,X_s)\\
-\sum_{i=1}^rT_v(\theta^1,\ldots,\theta^r,X_1,\ldots,\Cb_v(Z,X_i),\ldots,X_s),
\end{multline}
where $\theta^1,\ldots, \theta^r$ are one-forms and $Z,X_1,\ldots,X_r$ vector fields, both on $M$.
Moreover, let us denote by $\mathcal R$ the curvature tensor of $\nabla$ and $\hat R$ the curvature tensor of $\hat \nabla$, then
\begin{align}
\label{R1}{\mathcal R}_v(X^{\mathcal H},Y^{\mathcal H})Z=&\hat R_v(X,Y)Z+\Cb_v(i_v^{-1}({\mathcal V}[X^{\mathcal H},Y^{\mathcal H}]),Z)\\
\label{R2}{\mathcal R}_v(X^{\mathcal H},Y^{\mathcal V})Z=&-\hat P_v(X,Z,Y)+(\hat\nabla_X \Cb)_v(Y,Z)-\Cb_v(\hat P_v(X,v,Y),Z)\\
\label{R2b}{\mathcal R}_v(X^{\mathcal V},Y^{\mathcal H})Z=&\hat P_v(Y,Z,X)-(\hat\nabla_Y \Cb)_v(X,Z)+\Cb_v(\hat P_v(Y,v,X),Z)\\
\label{R3}{\mathcal R}_v(X^{\mathcal V},Y^{\mathcal V})Z=&\Cb_v(X,\Cb_v(Y,Z))-\Cb_v(Y,\Cb_v(X,Z))\nonumber\\
&+(\partial^\nu \Cb)_v(Y,Z,X)-(\partial^\nu \Cb)_v(X,Z,Y),
\end{align}
where $X,Y,Z$ are vector fields on $M$.
\end{proposition}
\begin{proof}
	The identity $\nabla_{Z^{\mathcal H}}T=\hat\nabla_ZT$ is checked straightforwardly from definitions. For \eqref{nablaZv}, observe that
	\[Z^{\mathcal V}(T(\theta^1,\ldots,\theta^r,X_1,\ldots,X_s))(v)=(\partial^\nu T)_v(\theta^1,\ldots,\theta^r,X_1,\ldots,X_s,Z)\]
	and for $i=1,\ldots,r$,
	\[\nabla_{Z^{\mathcal V}}\theta^i(X)=Z^{\mathcal V}(\theta^i(X))-\theta^i(\nabla_{Z^{\mathcal V}}X)=-\theta^i(\Cb(Z,X)).\]
	To check \eqref{R1}, observe that $d\pi([X^{\mathcal H},Y^{\mathcal H}])=[X,Y]$ and then ${\mathcal H}([X^{\mathcal H},Y^{\mathcal H}])=[X,Y]^{\mathcal H}$. Moreover, ${\mathcal V}([X^{\mathcal H},Y^{\mathcal H}])$ is tensorial in $X$ and $Y$. To obtain \eqref{R2}, take into account that
	\begin{multline}\label{first}
	(\nabla_{X^{\mathcal H}}(\nabla_{Y^{\mathcal V}}Z))(v)=\nabla_{X^{\mathcal H}}(\Cb(Y,Z))(v)=\hat\nabla^v_{X}(\Cb(Y,Z))\\=
(\hat\nabla_X\Cb)_v(Y,Z)+	\Cb_v(\hat\nabla_XY,Z)+\Cb_v(Y,\hat\nabla_XZ),
	\end{multline}
	where we have used the first identity of this proposition,
	\begin{equation}\label{second}
	(\nabla_{Y^{\mathcal V}}(\nabla_{X^{\mathcal H}}Z))(v)=\nabla_{Y^{\mathcal V}}(\hat\nabla_XZ)(v)=\hat P_v(X,Z,Y)+\Cb_v(Y,\hat\nabla_XZ),
	\end{equation}
	where we have used \eqref{nablaXVZ}. Finally, observe that $[X^{\mathcal H},Y^{\mathcal V}]$ is vertical, and 
	$[X^{\mathcal H},Y^{\mathcal V}](v)=i_v^{-1}(\hat\nabla^v_XY+\hat P_v(X,v,Y))$, which can be checked in coordinates using \eqref{Ngamma}. Then, using the last identity, one gets
	\begin{equation}\label{third}
	\nabla_{[X^{\mathcal H},Y^{\mathcal V}]}Z=\Cb(\hat\nabla_XY,Z)+\Cb(\hat P_v(X,v,Y),Z).
	\end{equation}
	Putting together \eqref{first}, \eqref{second} and \eqref{third}, one deduces \eqref{R2}. For \eqref{R2b}, just observe that ${\mathcal R}_v(X^{\mathcal V},Y^{\mathcal H})Z=-{\mathcal R}_v(Y^{\mathcal H},X^{\mathcal V})Z$ and apply \eqref{R2}. In order to check \eqref{R3}, just observe that using \eqref{nablaXVZ}, it follows
	\begin{align*}
	\nabla_{X^{\mathcal V}}(\nabla_{Y^{\mathcal V}}Z)(v)&=\nabla_{X^{\mathcal V}}(\Cb(Y,Z))(v)=(\partial^\nu\Cb)_v(Y,Z,X)+\Cb_v(X,\Cb(Y,Z)),\\
		\nabla_{Y^{\mathcal V}}(\nabla_{X^{\mathcal V}}Z)(v)&=(\partial^\nu\Cb)_v(X,Z,Y)+\Cb_v(Y,\Cb(X,Z)),
	\end{align*}
	and $[X^{\mathcal V},Y^{\mathcal V}]=0$.
	\end{proof}
\begin{corollary}\label{nablagvert}
	With the above notation, given a pseudo-Finsler metric $L:A\rightarrow \R$ and $g$ its fundamental tensor, one has that 
	\begin{equation}\label{nablag}
	\nabla_{X^{\mathcal V}}g(Y,Z)(v)=2 C_v(Y,Z,X)-g_v(\Cb(X,Y),Z)-g_v(Y,\Cb(X,Z)),
	\end{equation}
	where $X,Y,Z$ are vector fields on $M$. Moreover, if we fix the value of the anisotropic tensor $\nabla g$, and this anisotropic tensor is symmetric, then the tensor $\Cb$ is determined whenever it satisfies that $g_v(\Cb(X,Y),Z)$ is symmetric in $X,Y,Z$.
	\end{corollary}
\begin{proof}
	The identity \eqref{nablag} follows from \eqref{nablaZv} and using that $(\partial^\nu g)_v=2C_v$. The last claim is immediate.
	\end{proof}
It follows straightforwardly using the Christoffel symbols of the classical connections (see for example page 39 in \cite{BCS00}) that the Chern and  the Cartan linear connections determine the same anisotropic linear connection as do the Berwald and the Hashiguchi ones. For an elegant proof of this fact see \cite{Vit16}.
\begin{enumerate}
\item The Chern connection is determined by the Levi-Civita anisotropic connection $\nabla$ and $\Cb=0$,
\item The Berwald connection is determined by the Berwald anisotropic connection ($\nabla-\mathfrak{L}^\flat$) and $\Cb=0$,
\item The Cartan connection is determined by the Levi-Civita anisotropic connection $\nabla$ and $\Cb=C^\flat$
\item The Hashiguchi connection  is determined by the Berwald anisotropic connection ($\nabla-\mathfrak{L}^\flat$) and $\Cb=C^\flat$,
\end{enumerate}
where $C^\flat$ is determined by satisfying $g_v(C^\flat(X,Y),Z)=C_v(X,Y,Z)$ for all $X,Y,Z\in {\mathfrak X}(M)$.
\begin{corollary}
 Given a pseudo-Finsler metric $L:A\rightarrow \R$, let us consider only classical connections with the associated anisotro\-pic tensor $\Cb$ satisfying that the tensor $g_v(\Cb(X,Y),Z)$ is symmetric in $X,Y,Z$ and with torsion-free  associated anisotro\-pic connection $\hat\nabla$.	Then the following characterizations of the classical connections hold:
	\begin{enumerate}[(i)]
		\item the (linear) Chern connection is the only one such that \[\nabla_Xg(Y,Z)=2 C(Y,Z,i^{-1}({\mathcal V}X)),\]
			\item the (linear) Berwald connection is the only one such that \[\nabla_Xg(Y,Z)=2\mathfrak{L}(\pi_{TM}({\mathcal H}X),Y,Z)+2 C(Y,Z,i^{-1}({\mathcal V}X)),\]
				\item the  Cartan connection is the only one such that $\nabla_Xg(Y,Z)=0$,
					\item the Hashiguchi connection is the only one such that \[\nabla_Xg(Y,Z)=2\mathfrak{L}(\pi_{TM}({\mathcal H}X),Y,Z),\]
		\end{enumerate}
	in all the cases for every $X\in {\mathfrak X}(A)$ and $Y,Z\in {\mathfrak X} (M)$.
	\end{corollary}
\begin{proof}
	It follows from the above description of classical (linear) connections in terms of their associated anisotropic connection and anisotropic tensor $\Cb$, and the following facts:
	 $\hat \nabla g=0$ if $\nabla$ is the (anisotropic) Chern connection;
	 $\hat\nabla g=2\mathfrak{L}$ for the (anisotropic) Berwald connection (see (6.30) in \cite{Sh01});
	 when $\Cb=0$, $\nabla_{X^{\mathcal V}}g(Y,Z)(v)=2 C_v(Y,Z,X)$ and when $\Cb= C^\flat$, $\nabla_{X^{\mathcal V}}g(Y,Z)(v)=0$ 	(see \eqref{nablag}).
	\end{proof}
\begin{remark}
	Observe that with the definition of torsion $T$ for a classical linear connection given in \cite[Definition 2.6]{Vit16}, it holds that
	\begin{align*}
	T(X^{\mathcal H},Y^{\mathcal H})&={\mathcal T}(X,Y), &	T(X^{\mathcal V},Y^{\mathcal V})&=0.\\\
		T(X^{\mathcal H},Y^{\mathcal V})&=-{\mathfrak C}(Y,X),&
		T(X^{\mathcal V},Y^{\mathcal H})&={\mathfrak C}(X,Y).
	\end{align*}
	In order to prove this, observe that the operator $\mathscr{J}$ which appears in the definition of the torsion in \cite[Definition 2.6]{Vit16} has the following properties:  $\mathscr{J}(X^{\mathcal H})=X^{\mathcal V}$,  $\mathscr{J}(X^{\mathcal V})=0$, $\mathscr{J}[X^{\mathcal H},Y^{\mathcal H}]=[X,Y]^{\mathcal V}$, $\mathscr{J}[X^{\mathcal H},Y^{\mathcal V}]=0$ and $\mathscr{J}[X^{\mathcal V},Y^{\mathcal V}]=0$  (recall that $d\pi[X^{\mathcal H},Y^{\mathcal H}]=[X,Y]$,  $[X^{\mathcal H},Y^{\mathcal V}]$ is vertical and  $[X^{\mathcal V},Y^{\mathcal V}]=0$). The above formulas imply that a classical linear connection is torsion-free if and only if its associated anisotropic connection is torsion-free and $\mathfrak{C}=0$, which happens with Chern and Berwald connections, but not with Cartan and Hashiguchi ones.
\end{remark}
\section{Lie derivatives}\label{lieder}

Given  $Z\in \mathfrak{X}(M)$, the anisotropic Lie derivative ${\mathcal L}_Z$ is the tensor derivation obtained in Theorem \ref{existderiv} when you consider  the anisotropic derivation $\delta^v(X)=[Z,X]$, for every $v\in TM\setminus 0$ and $X\in \mathfrak{X}(M)$. Therefore, we get the anisotropic Lie derivative tensor ${\mathcal L}_Z T$ for every anisotropic tensor $T$. Let us see some particular cases. Let $L:A\rightarrow \R$ be a pseudo-Finsler metric. Then for every $X\in \mathfrak{X}(M)$, we have
\[(\mathcal L_XL)_V= X(L(V))-\partial^\nu L([X,V]),\]
where $V$ is any vector field that extends $v\in A$. Now as $\partial^\nu L(w)=2g_v(v,w)$,
 using the Chern connection and that $L(v)=g_v(v,v)$ we get
\[(\mathcal L_X L)_V= 2 g_V(\nabla^V_X V,V)-2g_V(V,[X,V])=2g_V(\nabla^V_VX,V).\]
It follows that $X$ is a Killing field of $L$ if and only if ${\mathcal L}_XL=0$ and conformal if and only if ${\mathcal L}_XL=f L$ for some function $f:M\rightarrow \R$ (see for example \cite[Proposition 6.1]{HJP15} and recall that the classical definition of a Killing (resp. conformal) field in Finsler geometry is a vector field such that its (possibly local) flow acts by isometries (resp. conformal maps) of the Finsler metric).
Moreover, 
\[(\mathcal L_X g)_v(Y,Z)= X(g_V(Y,Z))-g_V([X,Y],Z)-g_V(Y,[X,Z])-2 C_V(Y,Z,[X,V])\]
since $\partial^\nu g=2 C$. Using the Chern connection and \eqref{partialTexpres}, we get the tensorial expression
\[(\mathcal L_X g)_v(u,w)=g_v(\nabla^v_uX,w)+g_v(u,\nabla^v_wX)+2 C_v(u,w,\nabla^v_vX),\]
which gives another characterization of Killing fields:
\[g_v(\nabla^v_uX,w)+g_v(u,\nabla^v_wX)+2 C_v(u,w,\nabla^v_vX)=0,\]
for every $v\in A$ and $u,w\in T_{\pi(v)}M$.
Let us observe that given a diffeomorphism $\psi:M\rightarrow M$, we can define the pullback $\psi^*(T)$ of an anisotropic tensor $T\in  \mathfrak{T}^0_s(M,A)$ as the anisotropic tensor given by $\psi^*(T)_v(u_1,\ldots,u_s)=T_{\psi^*(v)}(\psi^*(u_1),\ldots,\psi^*(u_s))$, where $\psi^*$ is the differential of $\psi$ and $u_1,\ldots, u_s\in T_{\pi(v)}M$.
\begin{proposition}\label{Lieflow}
If $X\in \mathfrak{X}(M)$ and $T\in \mathfrak{T}^0_s(M,A)$, then 
\[\mathcal L_XT=\lim_{t\rightarrow 0}\frac{1}{t}(\psi^*_t(T)-T),\]
where $\psi_t$ is the (possibly local) flow of $X$.  
\end{proposition}
\begin{proof}
It follows the same steps as in the classical isotropic case (see for example \cite[Proposition 9.21]{Oneill} except for  the following:
\begin{multline*}
\lim_{t\rightarrow 0} \frac{1}{t}\left( T_{\psi^*_t(Y_p)}(V_{\psi_t(p)},W_{\psi_t(p)})-T_{Y_p}(V_p,W_p)\right)\\
=\lim_{t\rightarrow 0} \frac{1}{t}\left( T_{\psi^*_t(Y_p)}(V_{\psi_t(p)},W_{\psi_t(p)})-T_{Y(\psi_t(p))}(V_{\psi_t(p)},W_{\psi_t(p)})\right)\\
+\lim_{t\rightarrow 0} \frac{1}{t}\left( T_{Y(\psi_t(p))}(V_{\psi_t(p)},W_{\psi_t(p)})-T_Y(V_p,W_p)\right)\\
=\partial^\nu T_{Y_p}(V_p,W_p,\lim_{t\rightarrow 0} \frac{1}{t}(\psi^*_t(Y_p)-Y(\psi_t(p))))+X(T_Y(V,W))(p).
\end{multline*}
In the last equality, we have used that if a function $f:\Omega\subset \R^n\rightarrow \R$ is $C^1$ and $\phi,\varphi:(-\epsilon,\epsilon)\rightarrow \Omega$ are also $C^1$ with $\lim_{t\rightarrow 0}\phi(t)=\lim_{t\rightarrow 0}\varphi(t)=p$, then $\lim_{t\rightarrow 0}\frac{1}{t}(f(\phi(t))-f(\varphi(t))=df_p(\lim_{t\rightarrow 0}\frac{1}{t}(\phi(t)-\varphi(t)))$ (apply, for example, the mean value theorem). Finally observe that 
\[\lim_{t\rightarrow 0} \frac{1}{t}(\psi^*_t(Y_p)-Y(\psi_t(p)))=
\lim_{t\rightarrow 0} \psi^*_t(Y_p-\psi^*_{-t}(Y(\psi_t(p))))=-[X,Y]_p,\]
where we have used \cite[Proposition 1.58]{Oneill}.

\end{proof}
\section*{Acknowledgments}
This activity was supported  by
	the programme Young leaders in research 18942/JLI/13 by Fundaci\'on S\'eneca, Regional Agency for Science and Technology from the Region of Murcia and Spanish  MINECO/FEDER project reference
	MTM2015-65430-P


\begin{thebibliography}{10}
	
		\bibitem{Mat80}
	{\sc H.-H. Matthias}, {\em Zwei {V}erallgemeinerungen eines {S}atzes von
		{G}romoll und {M}eyer}, Bonner Mathematische Schriften [Bonn Mathematical
	Publications], 126, Universit\"at Bonn, Mathematisches Institut, Bonn, 1980.
	\newblock Dissertation, Rheinische Friedrich-Wilhelms-Universit{\"a}t, Bonn,
	1980.
	
		\bibitem{Jav14a}
	{\sc M.~A. Javaloyes}, {\em Chern connection of a pseudo-{F}insler metric as a
		family of affine connections}, Publ. Math. Debrecen, 84 (2014), pp.~29--43.
	
	\bibitem{Ja14b}
	{\sc M.~A. Javaloyes}, {\em Corrigendum to
		``{C}hern connection of a pseudo-{F}insler metric as a family of affine
		connections'' [mr3194771]}, Publ. Math. Debrecen, 85 (2014), pp.~481--487.
		
		\bibitem{JaSo14}
	{\sc M.~A. Javaloyes and B.~L. Soares}, {\em Geodesics and {J}acobi fields of
		pseudo-{F}insler manifolds}, Publ. Math. Debrecen, 87 (2015), pp.~57--78.	
	
		\bibitem{Sh01}
	{\sc Z.~Shen}, {\em Differential geometry of spray and {F}insler spaces},
	Kluwer Academic Publishers, Dordrecht, 2001.
	
	
		\bibitem{Yano57}
	{\sc K.~Yano}, {\em The theory of {L}ie derivatives and its applications},
	North-Holland Publishing Co., Amsterdam; P. Noordhoff Ltd., Groningen;
	Interscience Publishers Inc., New York, 1957.
	
		\bibitem{Lovas04}
	{\sc R.~L. Lovas}, {\em On the {K}illing vector fields of generalized metrics},
	SUT J. Math., 40 (2004), pp.~133--156.
	
	
		\bibitem{Oneill}
	{\sc B.~O'Neill}, {\em Semi-{R}iemannian geometry}, vol.~103 of Pure and
	Applied Mathematics, Academic Press, Inc. [Harcourt Brace Jovanovich,
	Publishers], New York, 1983.
	\newblock With applications to relativity.
		\bibitem{HJP15}
	{\sc J.~Herrera, M.~A. Javaloyes, and P.~Piccione}, {\em On a monodromy theorem
		for sheaves of local fields and applications}, Rev. R. Acad. Cienc. Exactas
	F\'{i}s. Nat. Ser. A Math. RACSAM, 111 (2017), pp.~999--1029.
	
	\bibitem{BCS00}
	{\sc D.~Bao, S.-S. Chern, and Z.~Shen}, {\em An introduction to
		{R}iemann-{F}insler geometry}, vol.~200 of Graduate Texts in Mathematics,
	Springer-Verlag, New York, 2000.
	
		\bibitem{Vit16}
	{\sc H.~Vit\'{o}rio}, {\em A unified approach to the theory of connections in
		{F}insler geometry}, Bull. Braz. Math. Soc. (N.S.), 48 (2017), pp.~317--333.
	

	

	


	
	
	
\end{thebibliography}
\end{document}